\documentclass[12pt]{amsart}

\newcommand{\be}{\begin{equation}}
\newcommand{\ee}{\end{equation}}
\newcommand{\bea}{\begin{eqnarray}}
\newcommand{\eea}{\end{eqnarray}}
\newcommand{\beaa}{\begin{eqnarray*}}
\newcommand{\eeaa}{\end{eqnarray*}}

\newcommand{\g}{{\bf g}}
\newcommand{\h}{{\bf h}}
\newcommand{\kk}{{\bf k}}

\advance\textheight by \topskip
\makeatletter

\@addtoreset{equation}{section}

\def\section{\@startsection {section}{1}{\z@}{-3.5ex plus -1ex minus
 -.2ex}{2.3ex plus .2ex}{\large\bf\centering}}
\def\subsection{\@startsection{subsection}{2}{\z@}{-3.25ex plus%
 -1ex minus -.2ex}{1.5ex plus .2ex}{\bf}}
\def\subsubsection{\@startsection{subsubsection}{3}{\z@}{-3.25ex plus%
 -1ex minus -.2ex}{1.5ex plus .2ex}{\sl}}
\makeatother
\begin{document}

 \baselineskip 17pt
\parindent 10pt
\parskip 9pt

\begin{flushright}
math.QA/0305285\\[3mm]
\end{flushright}

\vspace{1cm}
\begin{center}
{\Large {\bf Twisted Yangians and symmetric pairs}}\\ \vspace{1cm}
{\large N. J. MacKay}
\\
\vspace{3mm} {\em Department of Mathematics,\\ University of York,
\\York YO10 5DD, U.K.\footnote{email: {\tt nm15@york.ac.uk} }}
\end{center}

\begin{center}
{\em Talk presented at RAQIS'03, Annecy-le-Vieux, March 2003}
\end{center}

\vspace{0.2in} \centerline{\bf\small ABSTRACT} \centerline{
\parbox[t]{5in}{\small \noindent We describe recent work on the
twisted Yangians $Y(\g,\h)$ which arise as boundary remnants of
Yangians $Y(\g)$ in 1+1D integrable field theories, bringing out
the special role played by the requirement that $(\g,\h)$ be a
symmetric pair.}}

\vspace{0.2in}

\section{Introduction}

In a series of recent papers we examined the boundary principal
chiral model. The work began with the discovery that both the
classically integrable boundary conditions and the
vector-representation boundary $S$-matrices for the bulk $G$-model
were classified and parametrized by the symmetric spaces $G/H$
\cite{macka01}. These two points of view were linked by the
discovery \cite{deliu01} that the boundary remnant of the bulk
Yangian symmetry $Y(\g)$\footnote{In the PCM the symmetry is
$Y_L(\g)\times Y_R(\g)$, but in this talk we shall just look at
one copy of $Y(\g)$. Our conclusions should apply to any
integrable field theory with bulk $Y(\g)$ symmetry.} is the
twisted Yangian $Y(\g,\h)$ \cite{molev96,macka02}. We went on to
exploit our presentation of $Y(\g,\h)$ to give the spectral
decompositions of a variety of reflection or `$K$'-matrices (the
intertwiners of $Y(\g,\h)$-representations). In this talk we
summarize how $Y(\g)$ and $Y(\g,\h)$ fit together, emphasizing the
special role of the requirement that $(\g,\h)$ be a symmetric
pair. The focus is therefore on the algebra rather than the
physics -- an introduction to the latter may be found in
\cite{macsh}.

\vfill \pagebreak
\section{Yangians}

We begin with a brief recapitulation of the structure of $Y(\g)$,
largely drawn from \cite{macka02}.

Suppose our 1+1D integrable quantum field theory to have a compact
global symmetry group $G$, whose algebra $\g$ is generated by
conserved charges $Q_0^a$ with structure constants $f^a_{\;\;bc}$
and (trivial) coproduct $\Delta$,
\begin{equation}\label{cr0} \left[ Q_0^a , Q_0^b \right] = i
f^{a}_{\;\;bc} Q_0^c \hspace{0.3in}{\rm and}\hspace{0.3in}
\Delta(Q_0^a) = Q_0^a \otimes 1 + 1 \otimes Q_0^a \,.
\end{equation}
The Yangian \cite{drinf85} $Y(\g)$ is the larger symmetry algebra
generated by these and further non-local conserved charges
$Q_1^a$, where
\begin{equation}\label{cr1}
\left[ Q_0^a , Q_1^b \right] = i  f^{a}_{\;\;bc} Q_1^c
\hspace{0.3in}{\rm and}\hspace{0.3in} \Delta(Q_1^a) = Q_1^a
\otimes 1 + 1 \otimes Q_1^a + {1\over  2}f^{a}_{\;\;bc} Q_0^b
\otimes Q_0^c\,.
\end{equation}
The requirement that $\Delta$ be a homomorphism fixes\footnote{For
$\g\neq sl(2)$. For the general condition see
Drinfeld\cite{drinf85}.}
\begin{equation}\label{YSerre}
 f^{d[ab} [Q_1^{c]},Q_1^d] \; = \; {i\over {12}}
  \,f^{ap i} f^{bq j}f^{cr k}f^{ijk} \, Q_0^{(p} Q_0^q Q_0^{r)} \;,
\end{equation}
where $(\,)$ denotes symmetrization and $[\,]$ anti-symmetrization
on the enclosed indices, which have been raised and lowered freely
with the invariant metric $\gamma$.

The Yangian is a deformation of the polynomial algebra $\g[z]$:
with $Q_1^a=zQ_0^a$, the undeformed algebra would satisfy
(\ref{YSerre}) with the right-hand side zero -- that is, $z^2$
times the Jacobi identity. In $Y(\g)$, (\ref{YSerre}) acts as a
rigidity condition on the construction of higher $Q^a_n$ from the
$Q_1^a$.

There is an (`evaluation') automorphism
\begin{equation}\label{evalAM}L_\theta\,:\; Q_0^a \mapsto Q_0^a\,,\qquad Q^a_1 \mapsto
Q_1^a+ \theta { c_A \over 4i\pi} Q^a_0 \,,\end{equation} where
$c_A=C_2^\g(\g)$ is the value of the quadratic Casimir
$C_2^\g\equiv \gamma_{ab}Q_0^a Q_0^b$ in the adjoint
representation. (This normalization is chosen so that $\theta$ is
the particle rapidity, as we shall see later.) Thus any
representation $v$ of $Y(\g)$ may carry a parameter $\theta$: the
action of $Y(\g)$ on $v^\theta$ is that of $L_\theta(Y(\g))$ on
$v^0$. The $i$th fundamental representation $v_i^\theta$ of
$Y(\g)$ is in general reducible as a $\g$-representation, with one
of its irreducible components (that with the greatest highest
weight, where these are partially ordered using the simple roots)
being the $i$th fundamental representation $V_i$ of $\g$. In the
simplest cases (which include all $i$ for $\g=a_n$ and $c_n$),
$v_i^\theta=V_i$ as a $\g$-representation, and $Q_1^a=\theta { c_A
\over 4i\pi} Q^a_0$ upon it.

\vfill \pagebreak
\section{Twisted Yangians $Y(\g,\h)$}

A {\bf symmetric pair} $(\g,\h)$ is a (here compact and simple)
$\g$ together with a (maximal) subalgebra $\h\subset\g$ invariant
under an involutive automorphism $\sigma$. We shall write $\g=\h
\oplus \kk$, so that $\h$ and $\kk$ are the subspaces of $\g$ with
$\sigma$-eigenvalues $+1$ and $-1$ respectively, and $$
[\h,\h]\subset \h\,,\qquad [\h,\kk]\subset \kk\,, \qquad
[\kk,\kk]\subset\kk\,.$$ We shall use $a,b,c,...$ for general
$\g$-indices, $i,j,k,...$ for $\h$-indices and $p,q,r,...$ for
$\kk$-indices.

We define  \cite{deliu01,macka02} the {\bf twisted Yangian}
$Y(\g,\h)$ to be the subalgebra of $Y(\g)$ generated by
\bea\label{Q0q} &&Q_0^i\\{\rm and} && \widetilde{Q}_1^p \equiv
Q_1^p + {1\over 4} [C, Q_0^p] \label{Q1q}\,, \eea where $C\equiv
\gamma_{ij}Q_0^i Q_0^j$ is the quadratic Casimir operator of $\g$,
 restricted to $\h$.

$Y(\g,\h)$ is also a deformation, this time of the subalgebra of
(`twisted') polynomials in $\g[z]$ invariant under the combined
action of $\sigma$ and $z\mapsto -z$. Its defining feature is that
$Y(\g,\h)$ is a left\footnote{The analogous right co-ideal
subalgebra may be obtained by reversing the sign of the second
term in (\ref{Q1q}).} co-ideal subalgebra \cite{deliu01}, $\Delta
(Y(\g,\h))\subset Y(\g)\otimes Y(\g,\h)$. This is the crucial
feature which allows the boundary states to form representations
of $Y(\g,\h)$ while the bulk states form representations of
$Y(\g)$, and depends on the symmetric-pair property: \beaa
\Delta(\widetilde{Q}_1^p) & = & \Delta\left(Q_1^p +{1\over 4}
[C,Q_0^p]\right)\\ & = & Q_1^p\otimes 1 + 1\otimes Q_1^p + {1\over
4}[C,Q_0^p]\otimes 1 + 1 \otimes {1\over 4}[C,Q_0^p]\\&& + {1\over
2} f^p_{\;\;iq} Q_0^i\otimes Q_0^q +{1\over 2} f^p_{\;\;qi}
Q_0^q\otimes Q_0^i+{1\over 2}[\gamma_{ij}Q_0^i\otimes Q_0^j,
Q_0^p\otimes 1 + 1 \otimes Q_0^p]
\\ & = & \widetilde{Q}_1^p \otimes 1 + 1 \otimes
\widetilde{Q}_1^p + [\gamma_{ij}Q_0^i\otimes Q_0^j, Q_0^p\otimes
1]
\\ & = &\widetilde{Q}_1^p \otimes 1 + 1 \otimes \widetilde{Q}_1^p +
{1\over 2}[\Delta(C)-C\otimes 1 - 1 \otimes C, Q_0^p\otimes
1]\,.\\ & \subset & Y(\g)\otimes Y(\g,\h) \eeaa holds essentially
because for a symmetric pair the only non-zero structure constants
are $f^{ijk}$ and $f^{ipq}$, and fails for a general subalgebra
$\h\subset\g$.

Various particular cases of $(\g,\h)$ have been studied before
\cite{molev96}: for example, section 3.5 of the first paper of
\cite{molev96} describes $Y(gl(n),so(n))$.

\vfill \pagebreak
\section{Bulk and boundary scattering}

In both the bulk and the boundary case the (twisted) Yangian
symmetry allows one to determine the $S$-matrix up to an overall
scalar factor. In the boundary case, this will lead us to another
crucial implication of the symmetric-pair property.

In the bulk case, we first note that, in the scattering of two
particle multiplets $u^\phi$ and $v^\theta$, the asymptotic state
$u^\phi\otimes v^\theta$ is decomposable into a sum of
$\g$-representations on each of which the $S$-matrix acts as the
identity (because of the $\g$ symmetry). This state is in general
irreducible as a $Y(\g)$ representation, but it may become
reducible at certain special values of $\phi-\theta$ (at which the
$S$-matrix typically has a pole). In particular, the edge of the
physical strip lies at $\phi-\theta=i\pi$. The crossing symmetry
of the $S$-matrix requires that it project onto the scalar
representation of $\g$ at this value, while the $Q_1^a$ map the
scalar into the adjoint representation for general $\phi-\theta$.
The tensor product graph construction \cite{macka91,gould02}
(which we shall not recap here; see \cite{macka02} for a brief
introduction) makes it clear how the difference in the values of
$C_2^\g$ between these two representations fixes the pole value,
and ensures through the factor ${c_A \over 4i\pi}$ in
(\ref{evalAM}) that this pole is at $i\pi$.

The boundary $S$-matrix is determined as follows, in the simplest
cases for which $v^\theta$ is an irreducible $\g$-representation
(hereafter `irrep') $V$. Writing it (again up to an overall
factor, when it is usually known as the `reflection' or $K$-
matrix)\footnote{In certain cases, in which non-self-conjugate
$\g$-representations branch to self-conjugate
$\h$-representations, $v^\theta$ may be conjugated by $K$
\cite{macka01}.} as  $ K_v(\theta): v^\theta \rightarrow
v^{-\theta}$ and intertwining the $Q_0^i$ (that is, from the
physics point-of-view, requiring their conservation in boundary
scattering processes) requires that $$ K_v(\theta) Q_0^i = Q_0^i
K_v(\theta) $$ (in which by $Q_0^i$ we mean here its
representation on $V$) and thus that $K_v(\theta)$ act trivially
on $\h$-irreducible components of $V$. So we have $$K_v(\theta) =
\sum_{W \subset  V} \tau_W(\theta) P_W\,, $$ where the sum is over
$\h$-irreps $W$ into which $V$ branches, and $P_W$ is the
projector onto $W$.

To deduce relations among the $\tau_W$ we intertwine the
$\widetilde{Q}_1^p$. Recall that, on a $\g$-irreducible
$v^\theta$, the action of $Q_1^p$ is given by $ Q_1^p = \theta {
c_A \over 4i\pi} Q_0^p$, so that $$\langle W \vert\vert
K_v(\theta) \left( \theta { c_A \over 4i \pi} Q_0^p + {1\over
4}[C, Q_0^p] \right)\vert\vert W' \rangle = \langle W
\vert\vert\left( -\theta { c_A \over 4 i\pi}Q_0^p + {1\over 4}[C,
Q_0^p] \right) K_v(\theta)\vert\vert W' \rangle \,,$$ for
$W,W'\subset V$. Thus when the reduced matrix element $\langle W
\vert\vert Q_0^p\vert\vert W' \rangle\neq 0$ we have \be\label{BG}
{\tau_{W'}(\theta) \over \tau_{W}(\theta)} =
\left[\Delta\right]\,, \qquad{\rm where} \quad\left[ A \right]
\equiv {\frac{i\pi A}{c_A}+\theta \over \frac{i\pi A}{c_A}
-\theta}\ee and $\Delta = C(W)-C(W')$. To find the $W,W'$ for
which $\langle W \vert\vert Q_0^p\vert\vert W' \rangle\neq 0$ we
recall that $\kk$ forms an irrep $K$ of $\h$. A necessary
condition for (\ref{BG}) to apply is then that $W\subset K\otimes
W'$. Although not automatically sufficient, this is sufficiently
constraining in simple cases to enable us to deduce $K_v(\theta)$
\cite{deliu01,macka02}.

We can describe $K_v(\theta)$ by using a graph, in which the nodes
are the  $\h$-irreps $W$, linked by edges, directed from $W_i$ to
$W_j$ and labelled by $\Delta_{ij}$, whenever $W_i\subset K\otimes
W_j$. To calculate the labels, we first write $ C = \sum_i c_i
C_2^{\h_i}$, where $\h=\bigoplus_i \h_i$ is a sum of simple
factors $\h_i$ (and $C_2^{\h_i}$ is the quadratic Casimir of
$\h_i$). The point here is that $C$ was written in terms of
generators of $\g$: there will be non-trivial scaling factors
$c_i$, which may be computed by taking the trace of the adjoint
action of $C$ on $\g$ (where we fix $\gamma$ to be the identity
both on $\g$ and on each $\h_i$),  yielding $$ c_i= {c_A \over
C_2^{\h_i}(\h_i) + {{\rm dim}\, \kk\; \over {\rm dim}\, \h_i}
C_2^{\h_i}(\kk)}\,.$$ This has a highly non-trivial implication
for the boundary $S$-matrix. The analogue of the crossing relation
for bulk $S$-matrices is the `crossing-unitarity' relation
\cite{ghosh94}. One requirement of this is that, at the edge
$\theta=i\pi/2$ of the physical strip for the boundary $S$-matrix,
$K$ project onto the scalar representation of $\h$. In the graph
described above, $\langle W \vert\vert Q_0^p\vert\vert 1
\rangle\neq 0$ only for $W=K$, and so we must have $C(K)={1\over
2}c_A$, or
 $$ {C(K)\over c_A}={1\over c_A} \sum_i c_i C_2^{\h_i}(K)= \sum_i \left( {C_2^{\h_i}(\h_i) \over
C_2^{\h_i}(\kk)} + {{\rm dim}\, \kk\; \over {\rm dim}\,
\h_i}\right)^{-1}={1\over 2}\,.$$ That this holds, and does so
only for symmetric pairs, is a result of \cite{godd85} (also known
as the `symmetric space theorem' \cite{daboul96}) and we see once
again the centrality of this property.

\vfill\pagebreak
\section{Concluding remarks}

One reason why it is appropriate to emphasize the centrality of
the symmetric-pair property is that in a closely-related and
well-developed field, that of D-branes in group manifolds (a.k.a.\
the boundary WZW model), it appears not to be necessary; $\sigma$
may be {\em any} automorphism \cite{fffs99}. (Whether the triality
of $d_4$ enjoys some special status is unclear.)

A unifying principle for the exceptional algebras is the `magic
square', or equivalently the Cvitanovic-Deligne exceptional series
\cite{deligne96}, in terms of which both the boundary $K$-matrices
\cite{macka02} and bulk $R$-matrices \cite{westb02} have a unified
structure. Once again symmetric pairs are crucial in the
construction \cite{barto02}.

As promised, we have focused here on algebraic structures rather
than on physics, and so have not given examples of boundary
$S$-matrices or the spectra of boundary bound states which may be
deduced from them. Such calculations are generally very tough --
they involve complex fusion/bootstrap calculations -- and the
results are typically more opaque than in bulk cases
\cite{short03}. However, it is interesting to note that for the
classical Grassmannians $$ {SU(N)\over S(U(M)\times U(N-M))}\,,
\quad{SO(N)\over SO(M)\times SO(N-M)}\,,\quad {Sp(N)\over
Sp(M)\times Sp(N-M)} \;\;\;(N,M\;{\rm even})\,,$$ there seems to
be a set of boundary states with masses $$ m_a=m\sin{a\pi\over h}
\sin{(p-a)\pi\over h}\,,$$ where $m$ is a mass-scale and
$(p,h)=(M+1,N),\,(M,N-2),\,(M+2,N+2)$ for the three cases
respectively (taking, without loss of generality, $M\leq N/2$). In
the $N\rightarrow\infty$ limit, these are proportional to the
values of the quadratic Casimir operator in the $a$th fundamental
representations of (and so to the energy levels of the 0+1D
principal chiral model defined on) $SU(M),\,SO(M)$ or $Sp(M)$
respectively.

\vfill\pagebreak

\end{document}